\documentclass[11pt]{article}
\usepackage{amssymb}
\usepackage{amsfonts}
\usepackage{amsmath}
\usepackage{graphicx}

\newtheorem{theorem}{Theorem}

\newtheorem{corollary}{Corollary}

\newtheorem{lemma}{Lemma}

\newenvironment{proof}[1][Proof]{\noindent\textbf{#1.} }{\ \rule{0.5em}{0.5em} \medskip}

\def\mn{\medskip\noindent}
\def\bn{\bigskip\noindent}
\def\be{\begin{equation}}
\def\ee{\end{equation}}
\begin{document}

\title{Random Oxford Graphs} \author{by Jonah Blasiak\thanks{Work done in the summer of
2003 at a Cornell REU supported by the NSF.},\ Princeton U. \\ Rick
Durrett\thanks{Partially supported by NSF grants from the probability program (0202935)
and from a joint DMS/NIGMS initiative to support research in mathematical biology
(0201037).},\ Cornell U.} \maketitle

\maketitle
\begin{abstract}
Inspired by a concept in comparative genomics, we investigate properties of randomly chosen
members of $G_1(m,n,t)$, the set of bipartite graphs with $m$ left vertices, $n$ right
vertices, $t$ edges, and each vertex of degree at least one. We give asymptotic results for
the number of such graphs and the number of $(i,j)$ trees they contain. We compute
the thresholds for the emergence of a giant component and for the graph to be connected.
\end{abstract}

\section{Introduction}
Biologists use an {\it Oxford grid} to indicate the relationship between two genomes. It is a matrix with $g(i,j)=1$ if part of chromosome $i$ in the species $A$ is homologous to part of chromosome $j$ in species $B$. The corresponding Oxford graph is the bipartite graph obtained by letting the chromosomes of species $A$ be vertices on the left and chromosomes of species $B$ be vertices on the right and with an edge from $i$ on the left to $j$ on the right if $g(i,j)=1$. Figure 1 gives the Oxford graph
for the autosomes (non-sex chromosomes) of elephant and humans.

Let $G_1(m,n,t)$, the set of bipartite graphs with $m$ left vertices, $n$ right vertices,
$t$ edges, and each vertex of degree at least one. The graph in Figure 1 is a member of
$G_1(22,27,44)$ but is it a typical member of that set? To answer this question we will
examine properties of randomly chosen members of $G_1(m,n,t)$ and of related families of
bipartite graphs. We begin by asking how many such graphs there are. To answer this
question we will investigate the model $G^r(m,n,t)$: fix a vertex set $L$ of size $m$ and
$R$ of size $n$, and pick $t$ of the $mn$ edges between $L$ and $R$ with replacement
(picking the same edge multiple times is allowed).
As usual, we are interested in the behavior of these random graphs as $t$, $m$, and $n$
go to infinity; when using the symbols $\approx$, $\sim$, and $\to$ we are tacitly
assuming the results hold as $t$, $m$, and $n$ go to infinity.  Standard results for the
birthday problem (see e.g. page 83 of Durrett 1995) show that the probability no edge is
picked twice is $\approx \exp(-t^2/2mn)$, which converges to a positive limit if $t/m \to
\rho$ and $t/n \to \lambda$, so not much is changed by picking with replacement, except
that the next question becomes much easier to answer.

\mn Q. How big is $G^r_1(m,n,t)$, the subset of $G^r(m,n,t)$ in which each vertex has
degree at least one?

\mn
To relate this to the classical occupancy problem, consider an $m \times n$ array of boxes and throw in $t$ balls. Let $A$ be the event that each row has at least one ball and $B$ be the event that each column has at least one ball. It is easy to see that (thanks to sampling with replacement) the probability of $B$ is not affected by conditioning on the number of balls in each row, so $A$ and $B$ are independent. Using the multinomial distribution
$$
P(A) = \frac{1}{m^t} \sum_* \frac{ t! }{ i_1! \cdots i_m!}
$$
where the sum is over all $i_1, \ldots i_m\ge 1$ with $i_1 + \cdots + i_m = t$. To evaluate the sum we rewrite it as
$$
\frac{t! e^{am}} {m^t a^t} \sum_* \prod_{j=1}^m  e^{-a} \frac{ a^{i_j} }{ i_j!}
= \frac{t! e^{am}} {m^t a^t} P( Z_1 \ge 1, \ldots Z_m \ge 1, Z_1 + \cdots Z_m = t )
$$
where $Z_i$ are independent Poisson with mean $a$.

It is easy to see that $P( Z_1 \ge 1, \ldots Z_m \ge 1 ) = (1-e^{-a})^m$. $E(Z_i | Z_i \ge 1 ) = a/(1-e^{-a})$, so if we pick $a$ so that $a/(1-e^{-a}) = t/m$ and let $\sigma_a^2 = \hbox{var}\,(Z_i | Z_i \ge 1 )$ then
$$
P( Z_1 + \cdots Z_m = t | Z_1 \ge 1, \ldots Z_m \ge 1 ) \sim 1/\sqrt{2\pi \sigma_a^2 m}
$$
A similar analysis applies to $P(B)$ giving the following result.

\begin{theorem}
Let $a/(1-e^{-a}) = t/m$ and $b/(1-e^{-b}) = t/n$ and suppose that $t/m \to \lambda$, $t/n \to \rho$. Then
$$
|G^r_1(m,n,t)| = (nm)^t P(A) P(B) \sim (t!)^2
\frac {(e^a - 1)^m a^{-t} (e^b - 1)^n b^{-t}} { 2\pi \sigma_a \sigma_b \sqrt{mn} }
$$
\end{theorem}

\noindent
As a consequence of Theorem 1 and the birthday problem result we can calculate $|G_1(m,n,t)|$ up to a constant factor.

\begin{corollary}
Under the assumptions of Theorem 1,
$$
e^{-\rho\lambda} \le \liminf \frac{ |G_1(m,n,t)| }{ |G^r_1(m,n,t)| }
\le \limsup \frac{ |G_1(m,n,t)| }{ |G^r_1(m,n,t)| } \le 1
$$
\end{corollary}

Even more important than allowing us to count the graphs, the proof of Theorem 1 allows us
to relate our graphs to ones studied by Molloy and Reed (1995) and Newman, Strogatz, and
Watts (2001). Let $Y$ be a random variable with distribution given by
$$
P(Y = k) = \frac{1}{1-e^{-a}}\frac{e^{-a}a^k}{k!}, \quad\hbox{for $k \geq 1$}
$$
and $P(Y = k) = 0$ otherwise.  We will say $Y$ has a {\it truncated Poisson distribution}
with parameter $a$, or $\bar{\cal P}(a)$ for short.  This distribution is the
limiting degree distribution of a graph from $G^r_1(m,n,t)$ if parameter $a$ is
chosen correctly.  We choose $a$ by equating the means of the two distributions.  The
truncated Poission distribution has mean $a/(1-e^{-a})$ and the mean degree of a left
(right) vertex is $t/m$ ($t/n$).

We can now define a new graph model that mimics the degree distribution of vertices from
$G^r_1(m,n,t)$. Label the left vertices $l_1, l_2, \ldots ,l_m$ and the right vertices
$r_1, r_2, \ldots, r_n$.  Let $d(l_i)$, $i = 1, \ldots , m$ be independent $\bar{\cal
P}(a)$ random variables where $a/(1-e^{-a})=t/m$; let $d(r_i)$, $i = 1, \ldots , n$ be
independent $\bar{\cal P}(b)$ random variables where $b/(1-e^{-b}) = t/n$. Condition on
the sum of the $d(l_i)$ being $t$ and condition on the sum of the $d(r_i)$ being $t$.
Make a set $L'$ ($R'$) with $d(l_i)$ ($d(r_i)$) copies of vertex $l_i$ ($r_i$).  Pair up
elements in $L'$ with elements in $R'$ uniformly at random.  Finally, collapse the vertex
copies into a single vertex and let the vertex pairings determine the edges of the graph
(which may have multiple edges between vertices). Call the resulting random graph
$TP(m,n,t)$. It is clear that the $G^r_1(m,n,t)$ and $TP(m,n,t)$ random graph models have
the same degree distribution, and it is not surprising that models are, in fact, the
same.

\begin{lemma} \label{lemma1} The models $G^r_1(m,n,t)$ and $TP(m,n,t)$ are the same.
\end{lemma}

We give the proof in the appendix. To study the question of the existence of a giant
component in our graph, we begin with the general case in which the degrees of the $m$
left vertices have distribution $p_k$ and the degrees of the $n$ vertices on the right
have distribution $q_k$. If we examine the cluster of a given vertex $v$ on the left then
its first generation members (at distance one from $v$) will have distribution $p_k$, but
the number of children of a member of the first generation will not have distribution
$q_k$. A vertex on the right with degree $k$ is chosen in the first generation with
probability proportional to $kq_k$. If we let $\nu = \sum_k kq_k$ and $\bar q_k =
(k+1)q_{k+1}/\nu$ then the number of children of a child of $v$ will have distribution
$\bar q_k$ and mean $\bar\nu = \sum_k k \bar q_k$. Here we have shifted the distribution
by 1 to remove the edge that we arrived on (so that $v$ is not counted as its own
grandchild). Readers who are used to the Erd\"os-Renyi random graphs should note that if
$q_k$ is Poisson($\lambda$), then $\bar q_k$ is again Poisson($\lambda$).

Similar calculations apply to the third generation. The members of the second generation have size biased degree distributions $\bar p_k = (k+1)p_{k+1}/\mu$ where $\mu = \sum_k kp_k$ and this distribution has mean $\bar \mu$. As the reader can probably guess by analogy with branching processes,

\begin{lemma} The condition for the existence of a giant component is
$\bar \mu \cdot \bar \nu > 1$
\end{lemma}

Molloy and Reed (1995), who wrote the condition in the equivalent form $\sum_k k(k-2) p_k > 0$, proved this in the ordinary (unipartite case), essentially by showing that the branching process analogy gives an accurate approximation of cluster sizes.  Newman, Strogatz, and Watts (2001), motivated by studies of the structure of the world wide web, collaboration graphs of scientists, and Fortune 1000 company boards of directors, extended Molloy and Reed's results to directed and bipartite graphs. Since Newman, Strogatz, and Watts published in {\it Physical Review E}, they did not have to prove their results. Instead, like physicists, they wrote generating function equations that come from thinking of cluster formation as a branching process. As the reader can see from the description, Lemma 2 is almost a known result. Since we need some of the details in the proof of Theorem 4, we will give a detailed proof for the special case that appears in Theorem 2.

Our next step is to see what Lemma 2 says about our example. If $p_k$ is $\bar{\cal P}(a)$ then $\mu = a/(1-e^{-a})$ so
\be
\bar p_k = \frac{1-e^{-a}}{a} (k+1) e^{-a} \frac{a^{k+1}}{(k+1)!(1-e^{-a})} = e^{-a} \frac{a^{k}}{k!}
\ee
i.e., the Poisson distribution with mean $a$. A similar calculation shows $\bar q_k$ is the Poisson distribution with mean $b$, so the condition for the existence of a giant component is $ab>1$.

To compute the survival probability of the branching process, let $\phi_1$, $\phi_2$, $\psi_1$, and $\psi_2$ be the generating functions of $p_k$, $q_k$, $\bar p_k$, and $\bar q_k$ respectively. Consider our branching process, starting from one vertex on the left and conditioned on having one individual in the first generation. We call this the homogeneous branching process, because the different distribution at the first step has been eliminated. The number of offspring this individual has in the third
generation has generating function $\psi_2(\psi_1(z))$. To check the order of the composition note that if $N$ has distribution $\bar q_k$ ($N$ is the number of vertices in the second generation) and $X_1, X_2, \ldots$ are independent with distribution $\bar p_k$ ($X_1$ is the number of children of a second generation vertex) then
\be
E(z^{X_1 + \cdots X_N}) = \sum_{k=0}^\infty P(N=k) \psi_1(z)^k = \psi_2(\psi_1(z))
\ee
Let $\zeta_R$ be the smallest solution of $\psi_2(\psi_1(\zeta))=\zeta$ in $[0,1]$, i.e., the extinction probability of the homogeneous branching process. By considering the number of individuals in the first generation, it follows that the extinction probability for the branching process starting with one individual on the left is
$$
\xi_L = \sum_{k=1}^\infty p_k \, \zeta_R^k = \phi_1(\zeta_R)
$$
We define $\zeta_L$ and $\xi_R$ similarly.

\begin{theorem} Let $a/(1-e^{-a}) = t/m$ and $b/(1-e^{-b}) = t/n$ and suppose that $t/m \to \lambda$, $t/n \to \rho$. When $ab<1$ the largest cluster is $O(\log(m+n))$. A giant component appears when $ab>1$. The fraction of vertices it contains on the left and right are asymptotically $1-\xi_L$ and $1-\xi_R$. The second largest component is $O(\log(m+n))$.
\end{theorem}

To illustrate the phase transition we will consider some examples. In the human elephant comparison in Figure 1, $a=1.071$ and $b=1.593$ so $ab=1.707$. With a total of 49 vertices, it is hard to recognize a giant component, but there is one component with 13 human and 19 elephant vertices.  Figure 2 gives a comparison of human and colobine monkey, one of our  fairly close primate relatives, which has $a=0.503$, $b=0.605$ and $ab=0.305$. In agreement with subcritical designation, there are 12 components with 2 vertices, three with 3 vertices, one with 4, and one with 6. Figure 3 gives a comparison of the human and cat genomes that has $a=1.151$, $b=0.802$, and $ab=0.925$. Figure 4 compares humans and dogs, an example with $a=2.873$, $b=1.477$, and $ab=4.245$. The drastic difference in the graphs in Figures 3 and 4 is somewhat surprising since the evolutionary distance from humans to cats and dogs are the same. In the human-dog graph there is one giant component and three components of size 2. To lead into our next topic we ask: Does the number of small components in these random graphs agree with what we expect?

To get prepared for our next result, which will help us answer this question, we will give a second derivation of the threshold that is easy to believe but difficult to make rigorous. Suppose we are interested in some property of $G_1(m,n,t)$. Define $a$ and $b$ by $a/(1-e^{-a}) = t/m$ and $b/(1-e^{-b}) = t/n$. Let $G(M,N,p)$ be the random bipartite graph in which there are $M=t/a$ vertices on the left, $N=t/b$ on the right, and edges are independently chosen with probability $p = ab/t = a/N = b/M$.  $M$ and $N$ are defined this way so that after removing isolated vertices from each side we get a graph similar to one from $G_1(m,n,t)$.  The calculation is not difficult: the number of non-isolated vertices on the left, ${\cal M}$, has expected value
$$
E{\cal M} = M(1-( 1 - p)^N) \approx M(1-e^{-a}) = \frac{t}{a} (1-e^{-a}) = m,
$$
the number of non-isolated vertices on the right has $E{\cal N} = n$, and the number of edges, ${\cal E}$, has expected value $E{\cal E} = MNp = Nb = t$. Since all of the graphs in $G_1(m,n,t)$ have the same probability under $G(M,N,p)$.

\begin{lemma} The distribution of $G(M,N,p)$ conditioned on ${\cal M} = m$, ${\cal N}=n$, ${\cal E}=t$ is that of $G_1(m,n,t)$.
\end{lemma}

It is easy to show that when $t/m \to \rho$ and $t/n \to \lambda$, ${\cal M}$, ${\cal N}$, and ${\cal E}$, will with high probability differ from their expected values by $o(n)$. It is intuitively clear, but seems hard to show, that the vector $({\cal M},{\cal N},{\cal E})$ satisfies the local central limit theorem, so the conditioning ${\cal M}=m$, ${\cal N}=n$, ${\cal E}=t$ has probability $O(1/n^{3/2})$ and any property of $G(M,N,p)$ that has asymptotic probability $1-o(n^{-3/2})$ will be inherited by $G_1(m,n,t)$. Once one believes this, the threshold result follows easily. $G(M,N,P)$ has a giant component if
$$
1 < Mp \cdot Np = \frac{t}{a} \cdot \frac{t}{b} \cdot \left( \frac{ab}{t} \right)^2 = ab
$$

For a new example, consider the number of $(i,j)$ trees in the random graph, i.e., the number of trees with $i$ vertices on the left and $j$ vertices on the right. We let the tree size stay fixed while taking $m,n,t$ to infinity. Once one knows that the number of labeled bipartite $(i,j)$ trees is $i^{j-1} j^{i-1}$ (see e.g., Saltykov 1995), the expected number of $(i,j)$ trees in $G(M,N,p)$ can be derived by a calculation analogous to the standard one for trees in a unipartite random graph (see Bollob\'as (2001) Theorem 5.5).
$$
\frac{i^{j-1} j^{i-1}}{i! j!} \frac{(e^{-b} a)^j (e^{-a} b)^i} {p}
$$

Based on the reasoning above we expect that the corresponding result will hold for $G_1(m,n,t)$.

\begin{theorem} \label{expected tree thm}
In $G^r_1(m,n,t)$, the expected number of $(i,j)$ trees
$$
EA_{i,j} \to \frac{i^{j-1} j^{i-1}}{i! j!} \frac{(e^{-b} a)^j (e^{-a} b)^i t} {a b}
$$
\end{theorem}
Since the existence of $(i,j)$ trees on disjoint sets of vertices are asymptotically independent, we expect that the number of such trees will have asymptotically a Poisson distribution, but we have not tried to prove that.

To see what Theorem 3 says, we will consider our four previous examples and a comparison of the human and lemur genomes given in Figure 5, which is somewhat surprising since this example has $ab=1.771$ but no (1,1) or (2,1) trees. Table 1 compares the expected and observed number of (1,1), (2,1) and (1,2) trees. In general, there is good agreement between the observed and expected values. Two notable exceptions are the number of (1,1) trees in examples 4 and 5 where the expected values are 0.86 and 2.63 while the observed values are 3 and 0. If we assume that the number of trees has a Poisson distribution then the probability of three or more (1,1) trees in $G^r_1(22,38,67)$ is 0.097, while the probability of no (1,1) tree in $G^r_1(20,22,38)$ is 0.072.

Our final problem is to determine when the graph will be connected. For the Erd\"os-Renyi unipartite random graph $G(N,p)$ in which there are $N$ vertices and edges are independently present with probability $p$, the transition to connectivity occurs when $p \approx (\log N)/N$. To see this we note that the number of edges incident to vertex is asymptotically Poisson($Np$). If we let $p = c(\log N)/N$, the probability of an isolated vertex is $\approx 1/N^c$, so the expected value is large when $c<1$ and
goes to 0 if $c>1$. Isolated vertices prevent connectivity, so a second moment calculation shows that if $c<1$ the probability of connectivity goes to 0.

The result in the other direction is more difficult, since one must consider all of the ways in which the graph can fail to be connected. A simple calculation (see Bollob\'as 2001, p.~104) shows that if $p=\theta/N$ and $\theta = o(N^{1/2})$ then the expected number of trees with $v$ vertices, $T_v$, has
$$
E_p(T_v) \sim \frac{1}{\theta} \frac{v^{v-2}}{v!} (\theta e^{-\theta})^v
$$
From this we see that if $\theta = c\log N$ and $1/2 < c < 1$ then asymptotically there are isolated vertices, but no trees of size $v \ge 2$. Bollobas (2001), see Section 7.1, combines this estimate with the fact that the largest tree in a supercritical random graphs has $O(\log n)$ vertices to prove (see Theorem 7.3 on page 164) that if $\theta = \log N + x + o(1)$ then the probability $G(N,p)$ is connected approaches $\exp(-e^{-x})$.

Saltykov (1995) has considered a question closely related to the connectivity problem for the random bipartite graph $G(M,N,T)$ in which there are $M$ vertices on the left, $N$ vertices on the right, and $T$ edges. Suppose $M \ge N$. Let $\alpha = M/N$ and $\beta = (1-1/\alpha) \log N$. His main result asserts that if
$$
(1+1/\alpha)T = (M+N) \{ \log (M+N) + x + o(1) \}
$$
then the number of isolated vertices has asymptotically a Poisson distribution with mean
$$
\lambda = \frac{ e^{-x} (1+ e^{-\beta}) }{ 1+ 1/\alpha }
$$
Recalling $\alpha = M/N$, we see that the transition to connectedness occurs when $T \sim M \log (M+N)$.

The corresponding result for our bipartite random graphs is

\begin{theorem}
Define $c$ by $t = c\frac{mn}{m+n}log(m+n)$ and suppose $m/n \to \alpha$, a positive finite limit. The probability $G^r_1(m,n,t)$ is connected tends to 0 or 1 depending on whether  $c$ has a limit $<1$ or $>1$.
\end{theorem}
Note that our threshold is asymptotically $\frac{1}{1+\alpha} m \log (m+n)$. The difference in thresholds should not be surprising given the results for $E_p(T_v)$ cited above. Our threshold is for the disappearance of (1,1) trees rather than the absence of isolated vertices, so this occurs at a smaller value of $t$.

The remainder of the paper is devoted to proofs. We take the results in the same order as in the introduction.

\section{Proof of Corollary 1}

\noindent
{\bf Corollary 1} {\it Under the assumptions of Theorem 1,}
$$
e^{-\rho\lambda} \le \liminf \frac{ |G_1(m,n,t)| }{ |G^r_1(m,n,t)| }
\le \limsup \frac{ |G_1(m,n,t)| }{ |G^r_1(m,n,t)| } \le 1
$$

\begin{proof}
The inequality $|G_1(m,n,t)| \le |G^r_1(m,n,t)|$ is trivial and proves the result for $\limsup$. To prove the other result let $E = A \cap B$ be the event that there are no isolated vertices and let $F$ be the event that all edges chosen are distinct. Let $P$ denote probabilities under $G^r_1(m,n,t)$. From the thought experiment of sampling with replacement until we have $t$ distinct edges it is clear that $P(E|F) \ge P(E)$ because if a graph has no isolated vertices after the first $t$ edges are chosen, it will have no isolated vertices when $t$ distinct edges are chosen.  From this we get
$$
\frac{ |G_1(m,n,t)| }{ |G^r_1(m,n,t)|}
= \frac{P(E \cap F)}{P(E)} = \frac{ P(E|F) P(F) } {P(E)} \ge P(F)
$$
The result for $\liminf$ now follows from the result for the birthday problem cited in the introduction, which gives the limiting behavior of $P(F)$.
\end{proof}

\section{Proof of Theorem 2}

\noindent
{\bf Theorem 2} {\it Let $a/(1-e^{-a}) = t/m$ and $b/(1-e^{-b}) = t/n$ and suppose that $t/m \to \lambda$, $t/n \to \rho$. When $ab<1$ the largest cluster is $O(\log(m+n))$. A giant component appears when $ab>1$. The fraction of vertices it contains on the left and right are $1-\xi_L$ and $1-\xi_R$. The second largest component is $O(\log(m+n))$.}

\medskip
The first step is to make the connection between the cluster size and the total progeny in a branching process. To do this, we note that instead of making all of the choices in pairing the duplicated left and right vertices at once, we can do them sequentially. Suppose that we start with vertex $l_1$. We then choose $d(l_1)$ times without replacement from the duplicated set of right vertices $R'$. Let $f_1(r_j)$ be the number of times vertex $r_j$ is chosen and let ${\cal J}_1 = \{ j : f_1(r_j)>0\}$. For each $j \in {\cal J}_1$, choose $d(r_j) - f_1(r_j)$ times without replacement from the duplicated set of left vertices $L'$ minus the $d(l_1)$ copies of $l_1$. Let $f_2(l_j)$ be the number of times vertex $l_j$ is chosen, let ${\cal J}_2 = \{ j : f_2(l_j)>0\}$, etc. We continue this procedure until the cluster containing $l_1$ has been constructed. We then choose some vertex not in the cluster containing $l_1$, generate its cluster, and continue until the random graph has been constructed.

From the construction it should be clear that if $Y^m_k = |{\cal J}_k|$ is the number of vertices in generation $k$ (of a graph from $TP(m,n,t)$) then as $m\to\infty$, $\{ Y^m_k, k \ge 1 \}$ converges to the branching process described in the introduction. There are two differences between the growing cluster and the limiting branching process. The first is that the possible choices are dictated by the empirical sequence of degrees $d(l_1), \ldots d(l_m)$ and $d(r_1), \ldots d(r_n)$ rather than the truncated Poisson distributions. The second is that the set of available degrees changes as choices are made.

The first difference disappears as $m \to\infty$ since by the law of large numbers, the empirical distribution of degrees converges to the underlying theoretical distribution. To estimate the effect of the second, let $r_k$ be a probability distribution on the positive integers, let $\eta>0$, and let $W(\omega)$ be a nondecreasing function of $\omega \in (0,1)$ so that the Lebesgue measure $|\{ \omega: W(\omega)=k\}|=r_k$. We say that $W$ is the \emph{mass function} of distribution $r$. If we remove an amount of mass $\eta$ from the distribution and renormalize to get a probability distribution, then the result will be larger in distribution than $U = (W(\omega) | \omega < 1-\eta)$ and smaller in distribution than $V = (W(\omega) | \omega > \eta)$. Note that $EV \le EW / (1-\eta)$.

\medskip
{\it Subcritical Case.} Suppose $ab<1$. Pick $\eta>0$ so that $ab/(1-\eta) < 1$. Let $\hat p^m_k$ and $\hat q^m_k$ be the empirical distributions of the degrees of vertices on the left and on the right, let $\mu^*_m$ and $\nu^*_m$ be the means of these empirical distributions, and $\bar \mu^*_m = \sum_k (k-1) \hat p^m_k/\mu^*_m$ and $\bar \nu^*_m = \sum_k (k-1) \hat q^m_k/\nu^*_m$ be the means of the size biased distributions. Since $p_k$ and $q_k$ have finite second moments it follows from the law of
large numbers and (1) (pg. 4) that $\bar\mu^*_m \to a$ and $\bar\nu^*_m \to b$.

From the choice of $\eta$ it follows that if $m$ is large then until a fraction $\eta$ of vertices have been used up on either side, the growing cluster is dominated by a subcritical branching process. To estimate the growth of the cluster, we take the approach of Molloy and Reed (1995) and expose the cluster of right vertices one at a time, i.e., we pick one of the current set of active right vertices and go through two generations to identify the right vertices connected to it. The chosen right vertex is removed from the set of active vertices and the new ones are added; we call this a \emph{step}.  Vertices in early generations need not be exposed before vertices in later generations, as described at the beginning of the section; any active vertex may be exposed at each step.

To prove the lower bound on the critical value, we will show that if $ab<1$ then for large $m$ the largest cluster is $O(\log(m+n))$. Pick a right vertex at random and let $X$ be $-1$ plus the number of right vertices that can be reached in two steps in the branching process ($z^{X+1} = \psi_2(\psi_1(z))$ ). Assuming cluster growth is a branching process, this represents the change in the size of the set of active right vertices in one step of the construction. Let $S_\ell = S_0 + X_1 + \ldots + X_\ell$, where $X_i$ are independent with distribution $X$. When $S_0 =1$, $S_\ell$ gives the size of the active set of vertices after $\ell$ vertices in the cluster have been exposed. The random variable $\tau = \inf \{ \ell : S_\ell=0 \}$ has the same distribution as the total progeny of the homogeneous branching process starting from one right vertex.

In the limiting branching process $\kappa(\theta) \equiv Ee^{\theta X} < \infty$ for all $\theta$. Since $\kappa(0)=1$ and $\kappa'(0) = EX < 0$ in the subcritical case, there is a $\theta > 0$ so that $\kappa(\theta)<1$.
Therefore
\be
\label{kappabound}
P(\tau > k ) \le P( S_k \ge 1 ) \le E e^{\theta S_k} = \kappa(\theta)^k
\ee
so we have a bound on the total number of individuals in the branching process. To extend the last result to the growing cluster, we begin by observing that if $\hat X$ is the corresponding quantity for the empirical distribution then the strong law of large numbers implies $E\exp(\theta \hat X) \to E\exp(\theta X)$. If $X^\eta$ is the distribution that dominates choices made at any time before a fraction $\eta$ of the vertices have been used on the left or the right, then (from the discussions earlier) $E\exp(\theta X^\eta) \le E\exp(\theta \hat X)/(1-\eta)$.  So if $m$ is large and $\eta$ is small $E\exp(\theta X^\eta) < 1$. It follows from \ref{kappabound} that there is a $\gamma>0$ so that $P(\tau > k) \le e^{-\gamma k}$. If we take $k_0 = (2/\gamma) \log n$ then $P(\tau > k_0) \le n^{-2}$. This and the corresponding argument for left vertices proves that the largest cluster is $O(\log(m+n))$.

\medskip
{\it Supercritical Case.} Given distributions $\tilde d$ and $\bar d$, $\|\tilde d - \bar d\| = (1/2) \sum_k |\tilde d_k - \bar d_k|$ is the total variation distance. If $m$ is large and the fraction of vertices chosen on either side is at most $\eta$, then the cluster growth process dominates a branching process with offspring distributions $\tilde p_k$ and $\tilde q_k$ with $\| \tilde p - \bar p\| \le 2\eta$ and $\| \tilde q - \bar q\| \le 2\eta$ where $\bar p$ and $\bar q$ are the size biased degree distributions. Let $W_p$ be the mass function of $\bar p$. Among all distributions $\tilde p$ with $\|\tilde p - \bar p\| \le 2\eta$, the smallest one, $\bar p^\eta$, is the distribution with mass function $W_p^\eta$; $W_p^\eta(\omega) = W_p(\omega-2\eta), \; w \in (2\eta,1]$ and $W_p^\eta(\omega) = 0, \; w \in (0,2\eta]$. Define $W_q$, $W_q^\eta$, and $\bar q^\eta$ in the analogous way.

If we let $a_\eta$ and $b_\eta$ be the means of $\bar p^\eta$ and $\bar q^\eta$ then the dominated convergence theorem implies that as $\eta \to 0$, we have $a_\eta \to a$ and $b_\eta \to b$, so $a_\eta b_\eta > 1$ for small $\eta$. Now if $0 \le z\le 1$ we have
$$
\left| \sum_k \bar p^\eta_k z^k - \sum_k \bar p_k z^k \right| \le \sum_k |\bar p^\eta_k - \bar p_k | \le 4 \eta \to 0
$$
From this we see that if $\psi^\eta_1$ and $\psi^\eta_2$ are generating functions of $\bar p^\eta$ and $\bar q^\eta$ then, uniformly on $[0,1]$, we have $\psi^\eta_1 \to \psi_1$, $\psi^\eta_2 \to \psi_2$, and  $\psi^\eta_2(\psi^\eta_1) \to \psi_2(\psi_1)$. This uniform convergence implies that the smallest fixed point of $\psi^\eta_2(\psi^\eta_1)$ converges to that of $\psi_2(\psi_1)$, i.e., the extinction probability $\xi_R^\eta \to \xi_R$ as $\eta\to 0$. In a similar way we can conclude $\zeta_L^\eta \to \zeta_L$, $\xi_L^\eta \to \xi_L$, and $\zeta_R^\eta \to \zeta_R$.

To study the size of clusters, as in the previous proof, we expose them one right vertex at a time. When we expose the grandchildren of an active vertex, one of them might already be in the active set.  We call such an event a \emph{collision}.  If a collision occurs, instead of adding the grandchild to the active set (as is usually done), we remove it from the active set.  To show that this does not slow down the branching process too much, we must bound the number of collisions.  When we look at the left vertex children of a right vertex, we cannot encounter one we have seen before, because the first time a left vertex is visited, all of its other right vertex neighbors are added to the active set and all collisions are removed.
Note that $\bar p^\eta$ and $\bar q^\eta$ are concentrated on $\{ 0, \ldots, L \}$ where $L = \max \{ W_p^\eta(1), W_q^\eta(1) \}$. Thus until $\delta n$ vertices have been exposed on the right, the number of edges with an end in the active set is at most $\delta n L$.  The probability of picking one of these edges in the exposure of an active vertex is at most $\delta nL^2/(t - \delta nL^2) \equiv \gamma$.

Let $Z$ be the number of grandchildren in the branching process in which the first generation is according to $\bar q^\eta$ and the second according to $\bar p^\eta$. Let $Y$ be the distribution of grandchildren in the branching process modified to correct for collisions; $Y=Z - 2\cdot \text{Binomial}(\gamma,Z)$. Therefore if $\delta$ is small, $EY = a_\eta b_\eta (1-2\gamma) > 1$.

Let $X=Y-1$ and define $S_\ell$ as before. Since $EX>0$ the random walk has positive probability of not hitting 0, so there is positive probability that the cluster growth persists until there are at least $\delta m$ left vertices or $\delta n$ right vertices. To prove that we will get at least one such cluster with high probability, it is enough to show that with high probability all unsuccessful attempts will use up at most $O(\log(m+n))$ vertices. For this guarantees that we will get a large number of independent trails before using a fraction $\delta/2$ of vertices on either side.

The random variable $X$ is bounded so $\kappa(\theta)=Ee^{\theta X}<\infty$ for all $\theta$. $\kappa(\theta)$ is convex, continuous and has $\kappa'(0)=EX>0$, $\kappa(\theta) \sim P(X=-1) e^{-\theta} \to +\infty$ as $\theta\to-\infty$, so there is a unique $\lambda>0$ so that $\kappa(-\lambda)=1$. In this case $E\exp(-\lambda S_k)$ is a nonnegative martingale. Due to the possible removal of active vertices, the random walk may jump down by more than 1, but its jumps are bounded so the optional stopping theorem implies that the probability of reaching 0 from $S_0=x$ is $\le e^{-\lambda x}$.

The last estimate implies that the probability that the set of active vertices grows to size $(2/\lambda)\log n$ without generating a large cluster is $\le n^{-2}$. Routine  large deviations estimates for sums of independent random variables show that if $C$ is large, the probability that the sum of $C \log n$ independent copies of $X$ is $\le (2/\lambda)\log n$ is at most $n^{-2}$. Thus the probability of exposing more than $C \log n$ vertices and not generating a large cluster is $\le 2n^{-2}$. Combining this with the estimate for left clusters, we have our bound on unsuccessful attempts and can conclude that with high probability there is a large cluster.

To finish up now, let $\epsilon = \delta/L^2$. Since the maximum degree of any vertex is $L$, we can expose $\epsilon n$ right vertices without using up $\delta n$ vertices on either side. A routine large deviations estimate shows that
$$
P( S_{\epsilon n} \le (\epsilon n) EX/2 ) \le C e^{-cn}
$$
Consider now two vertices $i$ and $j$. If their clusters reach size $C \log n$ then the probability one of them will fail to continue until $\epsilon n$ right vertices have been exposed is $\le 4n^{-2}$. If the number of right vertices of their clusters reach size $\epsilon n$ and they have not already intersected, then with probability $ \ge 1 - 2C e^{-cn}$ each has an active set of size $\ge (\epsilon n) EX/2$.  The probability they will fail to intersect on the next step is exponentially small.  With probability tending to 1, all vertices in clusters larger than $C \log n$ belong to the giant component, and therefore the second largest component is $O(\log(m+n))$.

Our final task is to prove the claim about the fraction of vertices on the left and right that belong to the giant component.  Previous arguments have shown that if $\delta$ is small, the extinction probability for the comparison branching processes are $\approx \xi_L$. We have shown that membership in the giant component is essentially the same as belonging to a component of size $\ge C \log n$. Now, the probability of a collision before reaching size $C \log n$ is at most
\begin{equation}
(C \log n)^2 \cdot \frac{L^2}{n}
\label{collision}
\end{equation}
so if $1_{i \in G}$ is the indicator function that left vertex $i$ is part of a component of size $\ge C \log n$ then $E(1_{i \in G}) \approx 1-\xi_L$. When two clusters do not intersect, their growth is independent so (\ref{collision}) implies that
$$
\hbox{var\,}( \sum_{i=1}^m 1_{i \in G} ) \le m^2 (C \log n)^2 \cdot \frac{L^2}{n}
$$
Chebyshev's inequality implies
$$
\frac{1}{m} \left( \sum_{i=1}^m 1_{i \in G} - P(i \in G) \right) \to 0
$$
in probability and the desired result follows.

\section{Proof of Theorem 3}

\noindent
{\bf Theorem 3}
{\it In $G^r_1(m,n,t)$, the expected number of $(i,j)$ trees}
$$
EA_{i,j} \to \frac{i^{j-1} j^{i-1}}{i! j!} \frac{(e^{-b} a)^j (e^{-a} b)^i t} {a b}
$$

\begin{proof}
Let ${\cal T}$ be a fixed vertex labeled $(i,j)$ tree (left vertex labels are some subset of $\{1, 2, \ldots, m\}$ of size i), let $k=|E({\cal T})|=i+j-1$, and let $D$ be the event that it exists as a component of our random graph. Let $C(m,n,t)$ be the number of edge-labeled multigraphs belonging to $G^r_1(m,n,t)$.
$$
P(D) = \binom{t}{k} k! \frac{ C(m-i,n-j,t-k) }{ C(m,n,t) }
$$
The $\binom{t}{k} k!$ term comes from all the ways of labeling the edges of the
tree and dividing the labels between tree and non-tree edges.
From lemma \ref{lemma1}, we know
$$
C(m,n,t) = \sum^*_a \frac{t!}{ a_1! a_2! \ldots a_m!}
\sum^*_b \frac{t!}{b_1! b_2! \ldots b_n!}
$$
By symmetry it suffices to study the $m$ part of the equation. From the proof of Theorem 1, we have
$$
\sum^*_a \frac{t!}{ a_1! a_2! \ldots a_m!}  \approx  \frac{(e^{a}-1)^m t!}{a^t \sqrt{2 \pi \sigma^2_{a} m}}
$$
Thus $C(m-i,n-j,t-k)/C(m,n,t)$ is the product of two symmetric terms; the one containing $m$ is
\begin{equation} \label{long}
\frac{(e^{a'}-1)^{m-i} (t-k)!}{\sqrt{2 \pi \sigma^2_{a'} (m-i)} a'^{t-k} }
\div
\frac{(e^{a}-1)^{m} t!}{\sqrt{2 \pi \sigma^2_{a} m} a^{t} }
\end{equation}
where $a'$ is determined by $a'/(1-e^{-a'})=(t-k)/(m-i)$.

The expression above is equal to
\begin{equation*}
\left( \frac{e^{a'}-1}{e^a-1} \right)^{m-i}
\left( \frac{a}{a'} \right)^{t-k}
\cdot \frac{ \sigma_a }{ \sigma_{a'} }
\cdot \sqrt{ \frac{m}{m-i} }
\cdot \frac{a^k}{(e^a-1)^i}
\cdot \frac{(t-k)!}{t!}
\end{equation*}
Since $i$ and $k$ are fixed $a'$ tends to $a$ and $\sigma_{a'} \to \sigma_{a}$
\begin{equation}
\frac{ \sigma_a }{ \sigma_{a'} }
\cdot \sqrt{ \frac{m}{m-i} }
\cdot \left( \frac{e^{a'}-1}{e^{a}-1} \right)^{-i}
\cdot \left( \frac{a}{a'} \right)^{-k} \to 1
\label{eone}
\end{equation}
To complete the proof, we will show that
\begin{equation}
\left( \frac{e^{a'}-1}{e^{a}-1} \right)^m \left( \frac{a}{a'} \right)^t \to 1 \label{one}
\end{equation}
This enough since it implies
\begin{eqnarray} \label{expected tree eq}
P(D) & = & {t \choose k} k!
\frac{a^k (t-k)!}{(e^{a}-1)^{i} t!} \frac{b^k (t-k)!}{(e^{b}-1)^{j} t!} \\
& \sim & \frac {a^{j-1} b^{i-1}}{t^k} \frac{a^i}{(e^{a}-1)^{i}} \frac{b^j}{(e^{b}-1)^j}
\sim \frac{(e^{-b} a)^j (e^{-a} b)^i t} {m^i n^j a b} \nonumber
\end{eqnarray}
Multiplying this by $i^{j-1} j^{i-1} {m \choose i} {n \choose j}$, the number of vertex labeled $(i,j)$ trees on $(m,n)$ vertices, and taking limits gives Theorem 3.

To prove (\ref{one}) we use the definitions of $a$ and $a'$ to get
\be
\left(\frac{e^{a'}-1}{e^{a}-1} \right)^m \left( \frac{a}{a'} \right)^t =
\left(\frac{t}{m} \cdot \frac{m-i}{t-k} \cdot \frac{a'}{a} \right)^m e^{(a'-a)m} \left(\frac{a}{a'}\right)^{t}
\label{decomp}
\ee
To simplify these terms, we compute $a'-a$. Let $f(a) =  a/(1-e^{-a})$. The definition of the derivative implies
\begin{equation*}
a'-a  \sim \frac{f(a') - f(a)}{f'(a)} = \frac{1}{f'(a)} \left( \frac{t-k}{m-i} - \frac{t}{m} \right)
\end{equation*}
The next step is to note
\begin{equation}
\frac{t-k}{m-i} = \frac{t}{m}-\frac{k}{m}+ \frac{ti}{m^2} + O\left(\frac{1}{m^2}\right) \label{lina}
\end{equation}
and conclude that
\begin{equation}
a'-a \sim \frac{1}{f'(a)} \cdot \frac{\lambda i - k}{m}.
\label{deriv}
\end{equation}

Now the first term on the RHS of (\ref{decomp}) is
\begin{align}
& \left(\frac{m-i}{m}\right)^m  \cdot \left( 1 + \frac{k}{t-k} \right)^m \cdot
\left( 1 + \frac{a'-a}{a} \right)^m
\nonumber
\\
& = \left( 1 - \frac{i}{m} + \frac{k}{t-k} + \frac{a'-a}{a} + o(1/m) \right)^m
\label{firstterm}
\\
& \sim \left( 1 + \frac{1}{m}\left(-i+\frac{k}{\lambda}+\frac{\lambda i - k}{af'(a)}\right) \right)^m
\nonumber
\\
& \to \exp\left( (\lambda i - k)\left(-\frac{1}{\lambda}+\frac{1}{af'(a)}\right)\right)
\nonumber
\end{align}
if $t/m \to \lambda$.
By (\ref{deriv}) the second term on RHS of (\ref{decomp}) converges to
$\exp\left( (\lambda i - k)/f'(a)\right)$.
For the third term we write
$$
\left(\frac{a}{a'}\right)^{t} = \exp\left( -t \log\left( 1 + \frac{a'-a}{a} \right) \right)
$$
Using (\ref{deriv}) and expanding $\log(1+x) = x + O(x^2)$ shows that the third term converges to
$$
\exp\left(\frac{(\lambda i - k)(-\lambda)}{a f'(a)} + O(\frac{1}{m})\right)
\to \exp\left(\frac{(\lambda i - k)(-\lambda)}{a f'(a)}\right)
$$

Adding the three exponents gives
$$
(\lambda i - k) \left( - \frac{1}{\lambda} + \frac{1}{f'(a)} + \frac{(1-\lambda)}{a f'(a)} \right)
$$
We want to prove this is 0, so we can ignore the factor in front. Combining the fractions over a common
denominator, discarding that denominator, and recalling $\lambda = f(a)$ we have
$$
- af'(a) + (a+1-f(a))f(a)
$$
To check that this is zero, we note that differentiating $f(a) = a/(1-e^{-a})$ gives
\be
f'(a) = \frac{1}{(1-e^{-a})}  - \frac{ a e^{-a}}{(1-e^{-a})^2}
= \frac{f(a)}{a} + \frac{1}{a} \cdot (a-f(a))f(a)
\label{derivative}
\ee
and the proof is complete. \end{proof}

\section{Proof of Theorem 4}

\noindent
{\bf Theorem 4} {\it Define $c$ by $t = c\frac{mn}{m+n}log(m+n)$ and suppose $m/n \to \alpha$, a positive finite limit. The probability $G^r_1(m,n,t)$ is connected tends to 0 or 1 depending on whether $c$ has a limit $<1$ or $>1$.}

We can assume without loss of generality that $m \ge n$ and hence $\alpha \ge 1$. The first half of the proof is to establish:

\begin{lemma}
Under the assumptions of Theorem 4, if $c$ has a limit $< 1$ then the probability $G^r_1(m,n,t)$ is connected tends to 0.
\label{tendsto0}
\end{lemma}

\begin{proof}
Our first step is to show that the asymptotics in the previous section, which were derived under the assumption that $t$, $m$, and $n$ were all of the same order, continue to hold under the assumptions of Theorem 4.  To do this, it suffices to show that (\ref{eone}) and (\ref{one}) hold. We begin by noting that $t/m \to \infty$ implies $a\to\infty$ and $1-e^{-a} \to 1$, so $a \sim t/m$.  To verify (\ref{eone}) we observe that since $a\to\infty$, $\sigma_a^2/a \to 1$, and  $\sigma_a/\sigma_{a'} \to 1$.  In addition we will soon see that $a'-a \to 0$, and therefore $\left( \frac{e^{a'}-1}{e^{a}-1} \right)^{-i} \to e^{i(a-a')} \to 1$.

To prove (\ref{one}), we begin, as before, by computing $a'-a$.
As we have already noted
\be
a \sim \frac{t}{m} = \frac{cn}{m+n} \log(m+n) \to \frac{c}{1+\alpha} \log(m+n)
\label{asympa}
\ee
The fact that $a \sim t/m$ and the definition of $a$ implies that for large $m$
\be
\frac{t}{m} \ge a \ge \frac{t}{m} \left( 1 - e^{-t/2m} \right) \ge \frac{t}{m} \left( 1 - (m+n)^{-\epsilon} \right)
\label{esttm}
\ee
for some $\epsilon>0$. Since $a\to\infty$, we have $f'(a) \to 1$. Using this with (\ref{deriv}) and (\ref{lina}) it follows that
\be
a' - a \sim \frac{t-k}{m-i} - \frac{t}{m} = -\frac{k}{m}+ \frac{ti}{m^2} + O\left(\frac{1}{m^2}\right)
\label{diffa}
\ee
This leads to the asymptotic formula
\be
a' - a \sim \frac{c i}{1+\alpha} \frac{\log(m+n)}{m} - \frac{k}{m} \sim \frac{a i - k}{m}
\label{asympda}
\ee

Now we analyze the first term in the decomposition: using $t/m \sim a$ and (\ref{asympda}), (\ref{firstterm}) becomes
\be
\left( 1 + \frac{1}{m}\left(-i+\frac{k}{a}+\frac{a i - k}{af'(a)}\right) \right)^m
\to \exp\left( -i + \frac{i}{f'(a)} \right) \to 1
\label{firsttwo}
\ee

The second and third terms in (\ref{decomp}) are
$$
e^{(a'-a)m} \left(\frac{a}{a'}\right)^t = \exp\left( (a'-a)m -t \log \left( 1 + \frac{a'-a}{a} \right) \right)
$$
Expanding $\log(1+x) = x + O(x^2)$ the exponent becomes
$$
\left(m-\frac{t}{a}\right) (a'-a) -  t\; O \left( \frac{a'-a}{a} \right)^2
$$
(\ref{esttm}) implies that the absolute value of the first term is
$$
\le (m+n)^{-\epsilon} \cdot \frac{t}{a} |a'-a| \to 0
$$
by (\ref{diffa}) and $a \sim t/m$. To prove that the second term tends to 0, we note that $t/a \sim m$ and use (\ref{asympda}) and (\ref{asympa}). Thus we have
\be
e^{(a'-a)m} \left(\frac{a}{a'}\right)^t \to 1
\label{secondterm}
\ee
Combining this with (\ref{firsttwo}) gives (\ref{one}).

Let ${\cal T}_1,{\cal T}_2$ be fixed disjoint trees of size $(i, j)$.  Let $A_{i,j}$ be the number of $(i,j)$ trees that are components of our random graph, with $D_{{\cal T}}$ indicating whether ${\cal T}$ is a component. Writing $A_{i,j} = \sum_{\cal T} D_{\cal T}$, squaring and taking expected value we have
\begin{align}
E(A^2_{i,j}) & =   {m \choose i}{n \choose j} (i^{j-1}j^{i-1}) E(D_{{\cal T}_1})  \label{E(D)} \\
& + {m \choose i}{m-i \choose i}{n \choose j}{n-j \choose j} (i^{j-1}j^{i-1})^2 E(D_{ {\cal T}_1} D_{{\cal T}_2}). \nonumber
\end{align}
The last term counts the number of disjoint $(i,j)$ trees; overlapping trees contribute nothing to the sum.
To calculate $E(D_{{\cal T}_1} D_{{\cal T}_2})$, we note that calculations at the beginning of this section have shown
$$
\frac{ C(m-i,n-j,t-k) }{ C(m,n,t) } \sim \frac{a^k}{(e^{a}-1)^{i} t^k} \frac{b^k}{(e^{b}-1)^{j} t^k}
$$
so we have
\begin{align}
P(D_{{\cal T}_1}=1=D_{{\cal T}_2}) & = {t \choose 2k} (2k)! \frac{ C(m-2i,n-2j,t-2k) }{ C(m,n,t) } \nonumber \\
& \sim  t^{2k} \frac{a^{2k}}{(e^{a}-1)^{2i} t^{2k}} \frac{b^{2k}}{(e^{b}-1)^{2j} t^{2k}} \nonumber
\end{align}
where the ${t \choose 2k} (2k)!$ term comes from all the ways of labeling the edges of the
trees and dividing the labels between the two tree's edges and the other edges.  Recalling (\ref{expected tree eq}), we have that $E(D_{{\cal T}_1} D_{{\cal T}_2}) \sim E(D_{{\cal T}_1})^2$ and therefore (\ref{E(D)}) implies $E(A_{i,j}^2) \sim E(A_{i,j}) + E(A_{i,j})^2$.  We wish to show that $E(A_{1,1}) \to \infty$ so that we can conclude $E(A_{1,1}^2) \sim E(A_{1,1})^2$ and apply the second moment method.

To see this, observe that $a \leq t/m$.  Then the simplified expression for $E(A_{i,j})$ when $i=j=1$ is bounded as follows
\begin{align}
E(A_{1,1}) & \sim e^{-a} e^{-b} t  \geq e^{-t/n} e^{-t/m} t = (m+n)^{-c} t \nonumber \\
& = (m+n)^{-1-c} \cdot c m n \log(m+n) \label{ED}
\end{align}
Since $c < 1$ and $m/n \to \alpha$, a constant, this expression goes to infinity.
Now applying the second moment method yields $P(A_{1,1} = 0) \to 0$ which tells us that the probability of the existence of a $(1,1)$ tree goes to $1$,  and gives the desired result.
\end{proof}

Before tackling the other direction we need a preliminary result

\begin{lemma}
Let $Z$ have truncated Poisson distribution with mean $\lambda/(1-e^{-\lambda})$.
\begin{equation}
P( Z \le \lambda/2) \le \exp( -0.15 \lambda)
\label{expbd1}
\end{equation}
If $L>1/\ln 2$ then
\begin{equation}
P( Z \ge L \lambda ) \le \frac{1}{1-e^{-\lambda}}\exp( \lambda - L \lambda \ln 2)
\label{expbd2}
\end{equation}
\end{lemma}

\begin{proof}
Let $Z'$ be the Poisson distribution with mean $\lambda$.  The moment generating function is $E e^{\theta Z'} = \exp( \lambda(e^\theta - 1) )$, so if $\theta < 0$
$$
e^{\theta \lambda/2} P( Z' \le \lambda/2 ) \le \exp( \lambda (e^\theta - 1)).
$$
Taking $\theta = -\ln 2$
$$
P( Z' \le \lambda/2 ) \le \exp\left(  -\frac{\lambda}{2} (1 - \ln 2) \right)
$$
Since $\ln 2 \le 0.7$ and $P(Z \le \lambda/2 ) \le P( Z' \le \lambda/2 )$ the first result follows. For the second we note that if $\theta > 0$
$$
e^{\theta L\lambda} P( Z' \ge L \lambda ) \le \exp( \lambda (e^\theta - 1))
$$
Take $\theta=\ln 2$ and note that since $L>0$ we have
$$
P(Z \ge L \lambda ) = \frac{1}{1-e^{-\lambda}} P( Z' \ge L \lambda ) \le \frac{1}{1-e^{-\lambda}} \exp( \lambda - L \lambda \ln 2),
$$
the desired result.
\end{proof}

\begin{lemma}
Under the assumptions of Theorem 4, if $c$ has a limit $> 1$ then the probability $G^r_1(m,n,t)$ is connected tends to 1.
\end{lemma}

\begin{proof}
Under the assumptions of Theorem 4, $a \sim (cn/m+n) \log(m+n)$ and $b \sim (cm/m+n) \log(m+n)$. Let $r = \lim cn/(m+n)$ and $s = \lim cm/(m+n)$. Without loss of generality $s \ge r$, i.e., $n \ge m$. Our first step is to get an upper bound on the maximum degree of a vertex, $D$. By (\ref{expbd2}) with $\lambda = s \log(m+n)$
$$
P( D \ge L s \log(m+n) ) \le c \exp(s \log(m+n)(1 - L \ln 2)) = c (m+n)^{s- s L \ln 2}
$$
where $c =\frac{1}{1-(m+n)^{-s}} $. Taking $L = (2+s)/(s \ln 2)$ the right-hand side is $\le 2 (m+n)^{-2}$ for sufficiently large $m+n$. Assume for the rest of the proof that $D \le L \log(m+n)$.

The number of vertices in the first four generations is at most $N = \sum_{i=1}^4 (L \log (m+n))^i$. We will show that with high probability, $N$ is at least $O((log(m+n))^2)$ and this cluster will connect up to all others.  Using the trivial inequality $t \ge \max\{m,n\} \ge (m+n)/2$, the probability that two edges pick the same vertex in the first four generations (call this a collision, as before) is
$$
 \le N^2 \frac{D}{t} \le N^2 \cdot \frac{ 2L \log (m+n) }{m+n}
$$
This is too big to ignore but the probability of two or more collisions is
$$
\le N^4 \cdot \left( \frac{ 2L \log (m+n) }{m+n} \right)^2 \le C \frac{ (\log(m+n))^{18} }{ (m+n)^2 }
$$
so with high probability there is at most one collision in the first four generations of the cluster containing any vertex.

Our assumptions imply $r+s>1$, so we can pick $r' < r$ and $s' < s$ with $r' \le s'$ and $r' + s' \in (1,2)$. Pick $K$ so that $Kr'(0.15) > 2$. If $a \ge \ln 2$ (which will be true for large $m$), then in the associated branching process ($Z_i =$ the number of vertices in generation $i$)
$$
P( Z_1 \le K+1 ) = \frac{1}{1-e^{-a}} \sum_{k=1}^{K+1} e^{-a} \frac{a^k}{k!}
\le 2(K+1) e^{-a} a^K
$$
$a \sim r \log(m+n)$ so if $m$ is large
$$
P(Z_1 \le K+1) \le (n+m)^{-r'}
$$
By similar reasoning if $m$ is large
$$
P(\; Z_2 \le K+1 | Z_1 = j \; ) \le (n+m)^{-s'}
\quad\hbox{for all $j\ge 1$}
$$
From this it follows that
\begin{equation}
P( \max\{Z_1,Z_2\} \le K+1 ) \le (m+n)^{-(r'+s')}
\label{maxbd}
\end{equation}

So with high probability $Z_1$ or $Z_2$ is large and this implies $Z_4$ is large with high probability.  For $2 \le i \le 3$ divide individuals in generation $i$ into groups of size $K$. Since the sum of independent Poisson distributions is Poisson and the truncated Poisson distribution dominates the Poisson distribution, we may apply (\ref{expbd1}) to each group of size $K$.  
$$ P(\text{children of group} < K r' \log(m+n) / 2) $$
$$\le \exp(-0.15 K r' \log(m+n)) \le \frac{1}{(m+n)^2}$$

Trivially, the number of groups in generation $i$ is $\le (L \log (m+n))^i$ so

\begin{eqnarray*}
&&P( \; Z_{i+1} < \frac{Z_i}{K} \cdot \frac{K r' \log(n+m)}{2} \; | \; Z_i \ge K \; ) \le \frac{(L \log(m+n))^3}{(m+n)^2}
\end{eqnarray*}
Using this with (\ref{maxbd}) we can conclude that there is a constant $\delta>0$ for large $m$
$$
P(\; Z_4 < \delta (\log(m+n))^2 \; ) \le 4 (m+n)^{-(r'+s')}
$$
This shows that with high probability all clusters have size at least $\delta (\log(m+n))^2$. It follows from the proof of Theorem 2 that with high probability all clusters will grow to size $\delta n$ and connect. For readers who may be concerned with how the constants in that proof depend on $a$ and $b$ we note that all we need is a lower bound on the growth so for this phase of the argument, we can fix $a'<a$ and $b'<b$ with $a'b'>1$.  Theorem 2 does not apply when $a$ and $b$ are $O(log(m+n))$, but all we need is a lower bound, so it suffices to apply Theorem 2 with $a'$ and $b'$.
\end{proof}

\section{Appendix}

\begin{proof}
We will prove that the models $G^r_1(m,n,t)$ and $TP(m,n,t)$ are the same by looking at the distributions they induce on the set of edge labeled multigraphs. To do this, we will have to augment the model descriptions to label the edges. If we pick edges with replacement and label the edges in the order drawn then the set of outcomes $\Omega$, written as vectors of edges, has $(mn)^t$ elements and $G^r(m,n,t)$ is uniform over the subset $\Omega_0$ in which each vertex has degree at least one.

To label edges in $TP(m,n,t)$, first generate $L'$ and $R'$, the duplicated sets of vertices. Attach to the elements of $L'$ numbers chosen at random from $\{1, 2, ... t \}$ and call these edge-labels. Do the same independently for $R'$. Connect the element edge-labeled $i$ in $L'$ and the element edge-labeled $i$ in $R'$, and label this edge $i$.

Consider an outcome $w_0 \in \Omega_0$ with degrees $i_1, \ldots i_m$ and the left and $j_1, \ldots j_n$ on the right. By calculations in the introduction, the
probability that a graph in TP will have the same degrees as $w_0$ is
$$
\left. \frac{t!}{i_1! i_2! \cdots i_m!} \frac{t!}{j_1! j_2! \ldots j_m!} \right/ S(m,t) S(n,t)
$$
where $S(m,t)$ and $S(n,t)$ are normalizing constants that make the sum 1. Now $w_0$'s edge labels determine the edge labels incident to each vertex.  For each left vertex i, let $E_i^L$ be the set of edge labels incident to $i$ in $w_0$; similarly, let $E_j^R$ be the set of edge labels incident to right vertex $j$.  In order for TP to generate $w_0$, for each left vertex i, the labels of the set of vertices in $L'$ that collapse to i must be $E_i^L$ (but the order of the labels among the collapsing vertices doesn't matter).  A similar statement holds for the right vertices. The probability that vertices are labeled as described is
$$
\frac{i_1! i_2! \cdots i_m!}{t!} \frac{j_1! j_2! \ldots j_m!}{t!}
$$
so the edge labeled graphs generated by $TP(m,n,t)$ are also uniform on $\Omega_0$.
\end{proof}

\section*{References}

\mn
Bigoni, F. et al. (1997) Mapping homology between human and black and white colobine monkey chromosomes by flourescent in situ hybridization. {\it Am. J. Primatology.} {\bf 42}, 289--298

\mn
Bollob\'as, B. (2001) \textit{Random Graphs.} 2nd edition, Academic Press, New York

\mn
Breen, M. et al.~(1999)  Reciprocal chromosome painting reveals detailed regions of conserved synteny between the karyotypes of the domestic dog ({\it Canis familiaris}) and human. {\it Genomics} 61, 145-155

\mn
Newman, M. E. J., S. H. Strogatz, and D. J. Watts (2001) Random graphs with arbitrary degree
distributions and their applications. {\it Phys. Rev. E.} 64, article no. 026118.

\mn
Molloy, M. and B. Reed (1995) A critical point for random graphs with a given degree sequence.
{\it Random Structures Algorithms.} 6, 161–-179

\mn
Muller, S. et al.~(1999) Defining the ancestral karyotype of all primates by multidirectional chromosome painting between tree shrews, lemurs and humans. {\it Chromosoma.} 108, 393-400

\mn
Murphy, W.J. et al. (1999) Development of a feline whole genome radiation hybrid panel and comparative mapping of human chromosome 12 and 22 loci. {\it Genomics.} 57 (1999), 1-8

\mn
Saltykov, A. I. (1995) The number of components in a random bipartite graph. {\it Discrete Math. Appl.}
5, 515--523

\mn
Wienberg, J., et al. (1997) Conservation of human vs. feline genome organization revealed by reciprocal chromosome painting. {\it Cytogenetics and Cell Genetics.} 77, 211-217

\mn
Yang, F., et al. (2003) Reciprocal chromosome painting among human, aardvark, and elephant (superorder Afrotheria) reveals the likely eutherian ancestral karyotype. {\it Proc. Nat. Acad. Sci.} 100, 1062--1066

\clearpage

\begin{center}
\begin{picture}(285,290)
\put(20,270){$H_{17}$} \put(40,273){\line(1,0){20}}  \put(65,270){$E_{11}$}
\put(20,250){$H_{5}$}  \put(40,253){\line(1,0){20}}  \put(65,250){$E_{3}$}
\put(20,230){$H_{9}$}  \put(40,233){\line(1,0){20}}  \put(65,230){$E_{10}$}
\put(20,210){$H_{20}$} \put(40,213){\line(1,0){20}}  \put(65,210){$E_{23}$}
\put(20,170){$H_{10}$} \put(40,173){\line(1,0){20}}  \put(65,170){$E_{18}$}
\put(20,150){$H_{12}$} \put(40,153){\line(1,0){20}}  \put(65,150){$E_{4}$}
\put(20,130){$H_{22}$} \put(40,133){\line(1,0){20}}  \put(65,130){$E_{25}$}
\put(20,110){$H_{8}$}  \put(40,113){\line(1,0){20}}  \put(65,110){$E_{15}$}
\put(40,168){\line(2,-1){20}}
\put(40,148){\line(2,-1){20}}
\put(40,138){\line(2,1){20}}
\put(40,118){\line(2,1){20}}
\put(120,270){$H_{14}$} \put(140,273){\line(1,0){20}}  \put(165,270){$E_{9}$}
\put(140,258){\line(2,1){20}}
\put(120,250){$H_{15}$} \put(140,253){\line(1,0){20}} \put(165,250){$E_{17}$}
\put(140,238){\line(2,1){20}}
\put(120,230){$H_{4}$}  \put(140,233){\line(1,0){20}} \put(165,230){$E_{5}$}
\put(140,228){\line(2,-1){20}}
\put(120,210){$H_{18}$} \put(140,208){\line(2,-1){20}}
                         \put(165,210){$E_{20}$} \put(185,213){\line(1,0){30}}
\put(120,190){$H_{19}$} \put(140,193){\line(1,0){20}} \put(165,190){$E_{13}$}
\put(120,170){$H_6$} \put(140,168){\line(2,-1){20}}
\put(140,198){\line(2,1){20}} \put(140,188){\line(2,-1){20}}
                         \put(165,170){$E_{2}$}
\put(120,150){$H_{1}$}  \put(140,153){\line(1,0){20}} \put(165,150){$E_{19}$}
\put(140,158){\line(2,1){20}} \put(140,148){\line(2,-1){20}}
\put(120,130){$H_{21}$} \put(140,133){\line(1,0){20}} \put(165,130){$E_{21}$}
                         \put(165,110){$E_{1}$}
\put(120,90){$H_{3}$}   \put(140,93){\line(1,0){20}} \put(165,90){$E_{22}$}
\put(140,96){\line(2,1){20}} \put(140,90){\line(2,-1){20}}
\put(140,99){\line(2,3){20}} \put(140,87){\line(2,-3){20}}
                         \put(165,70){$E_{24}$}
                         \put(165,50){$E_{26}$}
\put(120,30){$H_{13}$}  \put(140,33){\line(1,0){20}} \put(165,30){$E_{16}$}
\put(140,38){\line(2,1){20}}
\put(220,230){$H_{7}$}  \put(240,233){\line(1,0){20}} \put(265,230){$E_{8}$}
\put(240,228){\line(2,-1){20}}
\put(220,210){$H_{16}$} \put(240,213){\line(1,0){20}} \put(265,210){$E_{12}$}
                        \put(265,190){$E_{27}$}
\put(220,170){$H_{2}$}  \put(240,173){\line(1,0){20}} \put(265,170){$E_{6}$}
\put(240,173){\line(1,0){20}} \put(240,176){\line(2,1){20}} \put(240,168){\line(2,-1){20}}
\put(240,179){\line(2,3){20}}
                         \put(265,150){$E_{14}$}
\put(220,130){$H_{11}$} \put(140,133){\line(1,0){20}} \put(265,130){$E_{7}$}
\put(240,133){\line(1,0){20}} \put(240,138){\line(2,1){20}}

\end{picture}
\end{center}

\mn
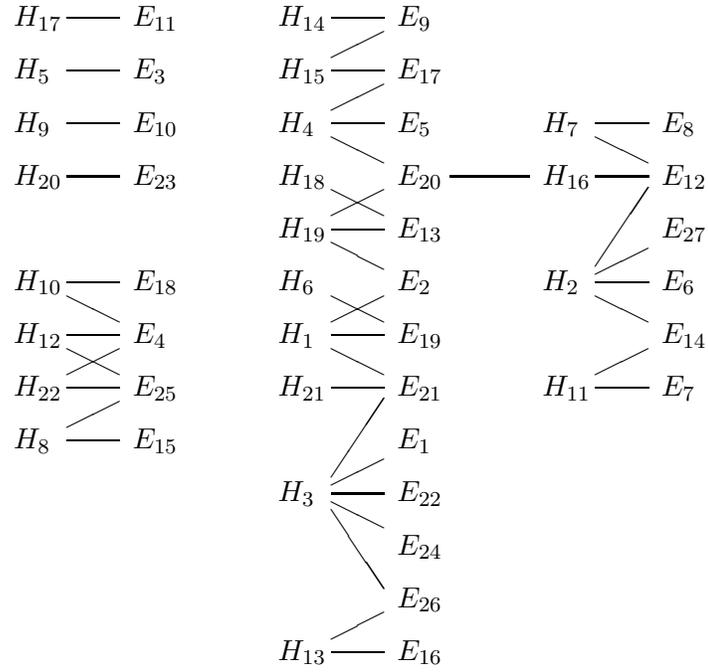
Figure 1. Comparison of elephant and human genomes.
Data from Yang et al. (2003). $m=22$, $n=27$, $t=44$, $a=1.126$, $b=1.654$, $ab=1.863$.
\clearpage

\begin{center}
\begin{picture}(310,250)
\put(20,190){$H_{4}$}
\put(20,170){$H_{5}$}
\put(20,150){$H_{6}$}
\put(20,130){$H_{7}$}
\put(20,110){$H_{8}$}
\put(20,90){$H_{9}$}
\put(20,70){$H_{11}$}
\put(20,50){$H_{12}$}
\put(20,30){$H_{13}$}
\put(65,190){$M_{1}$}
\put(65,170){$M_{2}$}
\put(65,150){$M_{3}$}
\put(65,130){$M_{4}$}
\put(65,110){$M_{9}$}
\put(65,90){$M_{10}$}
\put(65,70){$M_{14}$}
\put(65,50){$M_{7}$}
\put(65,30){$M_{19}$}
\put(40,193){\line(1,0){22}} \put(40,173){\line(1,0){22}} \put(40,153){\line(1,0){22}}
\put(40,133){\line(1,0){22}} \put(40,113){\line(1,0){22}} \put(40,93){\line(1,0){22}}
\put(40,73){\line(1,0){22}} \put(40,53){\line(1,0){22}} \put(40,33){\line(1,0){22}}
\put(120,190){$H_{16}$}
\put(120,170){$H_{18}$}
\put(120,150){$H_{20}$}
\put(165,190){$M_{15}$}
\put(165,170){$M_{20}$}
\put(165,150){$M_{21}$}
\put(140,193){\line(1,0){22}} \put(140,173){\line(1,0){22}} \put(140,153){\line(1,0){22}}
\put(120,110){$H_{14}$}
\put(120,90){$H_{15}$}
\put(165,100){$M_{6}$}
\put(142,113){\line(2,-1){20}} \put(142,93){\line(2,1){20}}
\put(120,50){$H_{21}$}
\put(120,30){$H_{22}$}
\put(165,40){$M_{16}$}
\put(142,53){\line(2,-1){20}} \put(142,33){\line(2,1){20}}
\put(265,190){$M_{17}$}
\put(265,170){$M_{20}$}
\put(220,180){$H_{16}$}
\put(242,183){\line(2,1){20}} \put(242,183){\line(2,-1){20}}
\put(220,130){$H_{3}$}
\put(220,110){$H_{19}$}
\put(265,130){$M_{12}$}
\put(265,110){$M_{18}$}
\put(240,113){\line(1,0){22}} \put(240,133){\line(1,0){22}}
\put(240,128){\line(2,-1){22}} \put(240,118){\line(2,1){22}}
\put(220,70){$H_{11}$}
\put(220,50){$H_{12}$}
\put(220,30){$H_{13}$}
\put(265,70){$M_{14}$}
\put(265,50){$M_{7}$}
\put(265,30){$M_{19}$}
\put(240,73){\line(1,0){22}} \put(240,51){\line(1,0){22}} \put(240,33){\line(1,0){22}}
\put(240,68){\line(2,-1){22}} \put(240,36){\line(2,1){22}}
\put(240,39){\line(2,3){22}}
\end{picture}
\end{center}

\mn
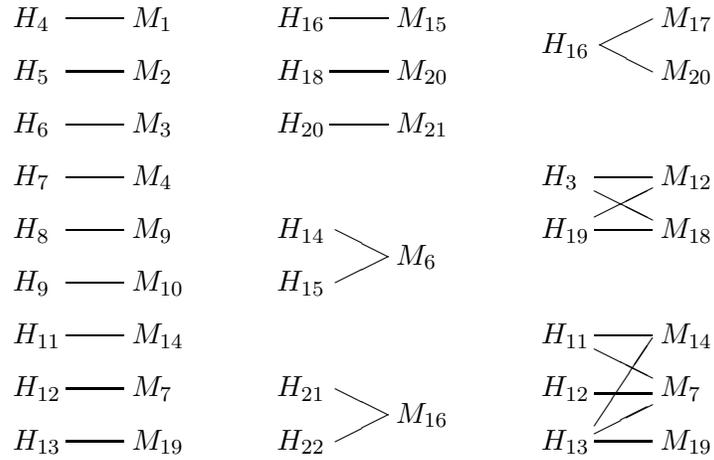
Figure 2. Comparison of human and colobine monkey ({\it Colobus guererza}) genomes. Data from Bigoni et al.~(1997). $m=22$, $n=21$, $t=28$, $a=0.581$, $b=0.685$, $ab = 0.397$.
\clearpage

\begin{center}
\begin{picture}(285,200)
\put(20,170){$H_6$} \put(40,173){\line(1,0){20}} \put(65,170){$B2$}
\put(20,150){$H_9$} \put(40,153){\line(1,0){20}} \put(65,150){$D4$}
\put(20,130){$H_{11}$} \put(40,133){\line(1,0){20}} \put(65,130){$D1$}
\put(20,110){$H_{17}$} \put(40,113){\line(1,0){20}} \put(65,110){$E1$}
\put(20,90){$H_5$} \put(40,93){\line(2,-1){20}}
\put(20,70){$H_{13}$} \put(40,73){\line(2,1){20}} \put(65,80){$A1$}
\put(20,50){$H_{14}$} \put(40,53){\line(2,-1){20}}
\put(20,30){$H_{15}$} \put(40,33){\line(2,1){20}} \put(65,40){$B3$}
\put(120,170){$H_{18}$} \put(140,168){\line(2,-1){20}}
\put(120,150){$H_{12}$} \put(140,153){\line(1,0){20}} \put(165,150){$D3$}
\put(140,148){\line(2,-1){20}} \put(140,138){\line(2,1){20}}
\put(120,130){$H_{22}$} \put(140,133){\line(1,0){20}} \put(165,130){$B4$}
\put(140,118){\line(2,1){20}}
\put(120,110){$H_{10}$} \put(140,113){\line(1,0){20}} \put(165,110){$D2$}
\put(120,70){$H_{20}$} \put(140,73){\line(1,0){20}} \put(165,70){$A3$}
\put(140,58){\line(2,1){20}}
\put(120,50){$H_2$} \put(140,53){\line(1,0){20}} \put(165,50){$C1$}
\put(140,38){\line(2,1){20}}
\put(120,30){$H_1$} \put(140,33){\line(1,0){20}} \put(165,30){$F1$}
\put(220,50){$H_4$} \put(240,53){\line(1,0){20}} \put(265,50){$B1$}
\put(240,38){\line(2,1){20}}
\put(220,30){$H_8$} \put(240,33){\line(1,0){20}} \put(265,30){$F2$}
\put(220,170){$H_{16}$} \put(240,173){\line(1,0){20}} \put(265,170){$E2$}
\put(240,168){\line(2,-1){20}} \put(240,158){\line(2,1){20}}
\put(220,150){$H_{19}$} \put(265,150){$E3$}
\put(240,148){\line(2,-1){20}} \put(240,138){\line(2,1){20}}
\put(220,130){$H_7$} \put(240,133){\line(1,0){20}} \put(265,130){$A2$}
\put(240,118){\line(2,1){20}}
\put(220,110){$H_3$} \put(240,113){\line(1,0){20}} \put(265,110){$C2$}
\put(220,90){$H_{21}$} \put(240,98){\line(2,1){20}}
\end{picture}
\end{center}

\mn
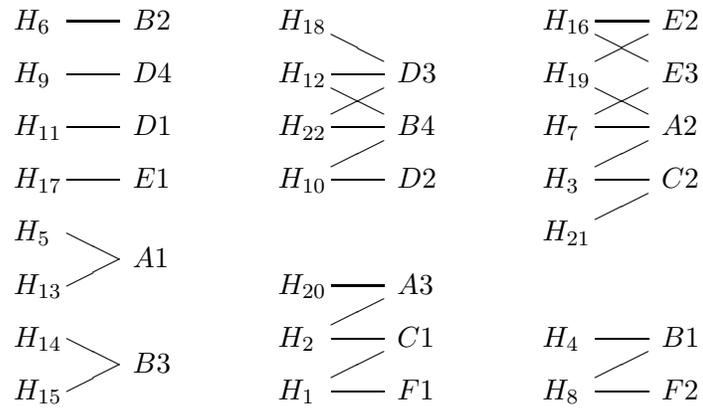
Figure 3. Comparison of human and cat genomes.
Data from Weinberg et al (1997) and Murphy et al (1999).
$m=22$, $n=19$, $t=32$, $a=1.151$, $b=0.802$, $ab = 0.925$
\clearpage

\begin{center}
\begin{picture}(310,280)

                                          \put(65,250){$D_{28}$}
                            \put(65,230){$D_4$}
                            \put(65,210){$D_2$}
\put(20,190){$H_{16}$} \put(40,188){\line(2,-1){22}} \put(40,198){\line(2,1){22}}
\put(40,183){\line(5,-6){22}}\put(65,190){$D_7$}
                        \put(65,170){$D_6$}
\put(20,150){$H_1$} \put(40,153){\line(1,0){22}} \put(65,150){$D_5$}
\put(40,148){\line(2,-1){22}} \put(40,156){\line(2,1){22}}  \put(40,163){\line(1,2){22}}
\put(40,143){\line(5,-6){22}}   \put(40,160){\line(5,6){22}}    \put(65,130){$D_{15}$}
                            \put(65,110){$D_{17}$}
                            \put(65,90){$D_{19}$}
\put(20,70){$H_2$} \put(40,73){\line(1,0){22}} \put(65,70){$D_{10}$}
\put(40,68){\line(2,-1){22}} \put(40,78){\line(2,1){22}}
\put(40,63){\line(5,-6){22}} \put(40,83){\line(5,6){22}}
                            \put(65,50){$D_{36}$}
                            \put(65,30){$D_{37}$}
\put(20,110){$D_{38}$} \put(25,142){\line(0,-1){20}}
\put(85,233){\line(1,0){22}} \put(110,230){$H_{10}$}
\put(85,213){\line(1,0){22}}\put(110,210){$H_{5}$}
\put(85,193){\line(1,0){22}} \put(110,190){$H_{18}$}
\put(110,170){$H_{17}$}
\put(110,150){$H_{4}$}
\put(110,130){$H_{7}$}
\put(110,110){$H_{11}$}
\put(85,73){\line(1,0){22}} \put(110,70){$H_{12}$}
\put(110,50){$H_{22}$}
\put(85,218){\line(2,1){22}} \put(85,158){\line(2,1){22}} \put(85,138){\line(2,1){22}}
\put(85,68){\line(2,-1){22}} \put(85,228){\line(2,-1){22}} \put(85,248){\line(2,-1){22}}
\put(85,100){\line(1,2){22}}
\put(82,170){\line(4,-5){25}} \put(82,150){\line(4,-5){25}}
\put(84,125){\line(1,-2){23}}
\put(155,270){$D_{12}$}
\put(155,250){$D_{35}$}
\put(155,230){$D_{11}$}
\put(155,210){$D_{1}$}
\put(155,190){$D_{9}$}
\put(155,170){$D_{3}$}
\put(155,150){$D_{13}$}
\put(155,130){$D_{14}$}
\put(155,110){$D_{16}$}
\put(155,90){$D_{18}$}
\put(155,70){$D_{21}$}
\put(155,50){$D_{27}$}
\put(155,30){$D_{26}$}
\put(130,218){\line(2,1){22}} \put(130,198){\line(2,1){22}} \put(130,178){\line(2,1){22}} \put(130,158){\line(2,1){22}}
\put(130,128){\line(2,-1){22}} \put(130,108){\line(2,-1){22}} \put(130,70){\line(2,-1){22}} \put(130,48){\line(2,-1){22}}
\put(130,153){\line(1,0){22}} \put(130,133){\line(1,0){22}} \put(130,113){\line(1,0){22}}
\put(130,103){\line(5,-6){22}} \put(130,66){\line(5,-6){20}} \put(127,208){\line(4,-5){25}}
\put(130,148){\line(5,-1){67}}
\put(200,130){$D_{32}$}
\put(200,250){$H_{6}$}
\put(200,230){$H_{3}$}
\put(200,210){$H_{9}$}
\put(200,190){$H_{19}$}
\put(200,170){$H_{15}$}
\put(200,150){$H_{8}$}
\put(200,110){$H_{13}$}
\put(200,70){$H_{14}$} \put(220,73){\line(1,0){22}} \put(245,70){$D_{8}$}
\put(200,50){$H_{20}$} \put(220,53){\line(1,0){22}} \put(245,50){$D_{24}$}
\put(200,30){$H_{21}$} \put(220,33){\line(1,0){22}} \put(245,30){$D_{31}$}
\put(175,253){\line(1,0){22}} \put(175,213){\line(1,0){22}} \put(175,173){\line(1,0){22}} \put(175,153){\line(1,0){22}}
\put(175,268){\line(2,-1){22}} \put(175,228){\line(2,-1){22}} \put(175,208){\line(2,-1){22}}
\put(175,218){\line(5,6){22}} \put(175,118){\line(5,6){22}}
\put(175,198){\line(2,1){22}}
\put(245,270){$D_{34}$}
\put(220,243){\line(5,6){22}}
\put(245,250){$D_{23}$}
\put(245,230){$D_{33}$}
\put(245,210){$D_{20}$}
\put(245,170){$D_{30}$}
\put(245,150){$D_{29}$}
\put(245,130){$D_{25}$}
\put(245,110){$D_{22}$}
\put(220,233){\line(1,0){22}} \put(220,173){\line(1,0){22}} \put(220,153){\line(1,0){22}} \put(220,113){\line(1,0){22}}
\put(220,228){\line(2,-1){22}} \put(220,148){\line(2,-1){22}}
\put(220,238){\line(2,1){22}} \put(220,198){\line(2,1){22}} \put(220,118){\line(2,1){22}}
\end{picture}
\end{center}

\mn
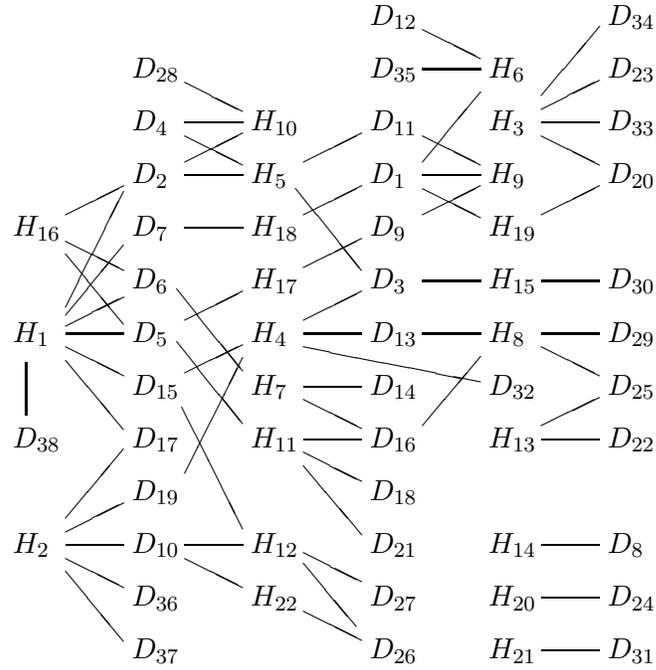
Figure 4. Comparison of the human and dog genomes. Data from Breen et al. (1999).
$m=22$, $n=38$, $t=67$, $a=2.873$, $b=1.477$, $ab=4.245$.
\clearpage

\begin{center}
\begin{picture}(310,250)
\put(20,190){$L_{1}$}
\put(20,150){$L_{9}$}
\put(20,130){$L_{14}$}
\put(20,110){$L_{18}$}
\put(20,70){$L_{17}$}
\put(20,50){$L_{5}$}
\put(20,30){$L_{13}$}
\put(65,210){$H_{9}$}
\put(65,190){$H_{3}$}
\put(65,170){$H_{21}$}
\put(65,130){$H_{1}$}
\put(65,70){$H_{10}$}
\put(65,50){$H_{12}$}
\put(65,30){$H_{22}$}
\put(40,193){\line(1,0){22}} \put(40,133){\line(1,0){22}} \put(40,73){\line(1,0){22}}
\put(40,53){\line(1,0){22}} \put(40,33){\line(1,0){22}}
\put(40,198){\line(2,1){22}} \put(40,118){\line(2,1){22}}\put(40,38){\line(2,1){22}} \put(40,58){\line(2,1){22}}
\put(40,188){\line(2,-1){22}} \put(40,148){\line(2,-1){22}} \put(40,48){\line(2,-1){22}}
\put(120,150){$L_{1}$}
\put(120,130){$L_{20}$}
\put(120,110){$L_{12}$}
\put(120,70){$L_{7}$}
\put(120,30){$L_{10}$}
\put(165,170){$H_{18}$}
\put(165,150){$H_{6}$}
\put(165,130){$H_{4}$}
\put(165,110){$H_{19}$}
\put(165,90){$H_{15}$}
\put(165,70){$H_{14}$}
\put(165,50){$H_{16}$}
\put(165,30){$H_{17}$}
\put(140,153){\line(1,0){22}} \put(140,133){\line(1,0){22}} \put(140,113){\line(1,0){22}}
\put(140,73){\line(1,0){22}} \put(140,33){\line(1,0){22}}
\put(140,158){\line(2,1){22}} \put(140,118){\line(2,1){22}} \put(140,78){\line(2,1){22}} \put(140,38){\line(2,1){22}}
\put(140,148){\line(2,-1){22}} \put(140,68){\line(2,-1){22}}
\put(140,83){\line(5,6){22}}
\put(210,150){$L_{11}$}
\put(210,130){$L_{6}$}
\put(210,90){$L_{2}$}
\put(210,70){$L_{16}$}
\put(185,153){\line(1,0){22}} \put(185,133){\line(1,0){22}} \put(185,93){\line(1,0){22}}
\put(255,130){$H_{20}$}
\put(255,110){$H_{2}$}
\put(255,90){$H_{8}$}
\put(230,133){\line(1,0){22}} \put(230,93){\line(1,0){22}}  \put(230,78){\line(2,1){22}} \put(230,98){\line(2,1){22}}
\put(230,128){\line(2,-1){22}}
\put(210,230){$L_{19}$}
\put(210,210){$L_{3}$}
\put(210,190){$L_{15}$}
\put(255,230){$H_{7}$}
\put(255,210){$H_{5}$}
\put(230,233){\line(1,0){22}} \put(230,213){\line(1,0){22}}
\put(230,218){\line(2,1){22}} \put(230,198){\line(2,1){22}}
\put(120,210){$L_{4}$}
\put(165,230){$H_{13}$}
\put(165,210){$H_{11}$}
\put(140,218){\line(2,1){22}}
\put(140,213){\line(1,0){22}}
\end{picture}
\end{center}

\mn
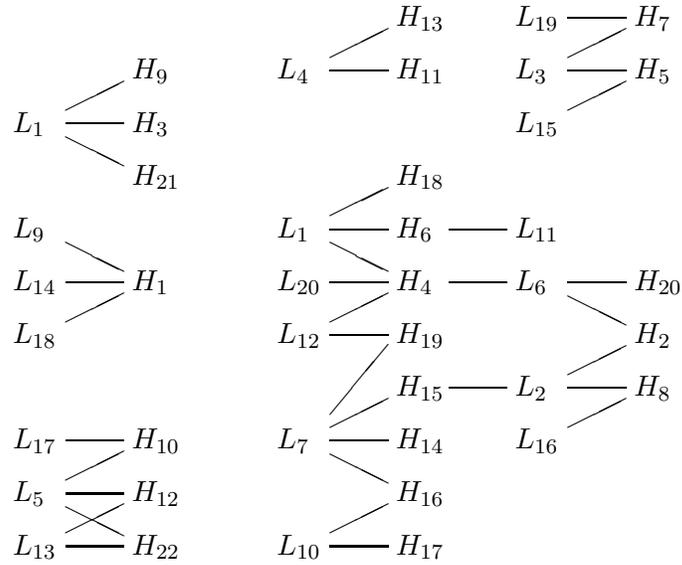
Figure 5. Comparison of lemur ({\it Eulemur macao macao}) and human genomes. Data from M\"uller et al.~(1997). $m=20$, $n=22$, $t=38$, $a=1.458$, $b=1.214$, $ab=1.771$.
\clearpage

\begin{center}
\begin{tabular}{cccccc}
Example & 1 & 2 & 3 & 4 & 5 \cr
human & elephant & monkey & cat & dog & lemur \cr
$ab$       & 1.71 & 0.30 & 0.93 & 4.25 & 1.77 \cr
$EA_{1,1}$ & 3.06 & 9.23 & 4.53 & 0.86 & 2.63 \cr
obs        & 4    & 12 & 4 & 3 & 0 \cr
$EA_{2,1}$ & 0.33 & 1.69 & 1.17 & 0.04 & 0.37 \cr
obs        & 0    & 2 & 2 & 0 & 0 \cr
$EA_{1,2}$ & 0.83 & 1.26 & 0.57 & 0.28 & 0.57 \cr
obs        & 0    & 1 & 0 & 0 & 1
\end{tabular}
\end{center}

\bn
Table 1. Expected number of trees of various sizes compared with the number observed in our five examples.

\end{document}